# ON THE $\mathbb{L}_P$-ERROR OF MONOTONICITY CONSTRAINED ESTIMATORS

By Cécile Durot

*Université Paris Sud*

We aim at estimating a function $\lambda:[0,1] \to \mathbb{R}$, subject to the constraint that it is decreasing (or increasing). We provide a unified approach for studying the $\mathbb{L}_p$-loss of an estimator defined as the slope of a concave (or convex) approximation of an estimator of a primitive of $\lambda$, based on $n$ observations. Our main task is to prove that the $\mathbb{L}_p$-loss is asymptotically Gaussian with explicit (though unknown) asymptotic mean and variance. We also prove that the local $\mathbb{L}_p$-risk at a fixed point and the global $\mathbb{L}_p$-risk are of order $n^{-p/3}$. Applying the results to the density and regression models, we recover and generalize known results about Grenander and Brunk estimators. Also, we obtain new results for the Huang–Wellner estimator of a monotone failure rate in the random censorship model, and for an estimator of the monotone intensity function of an inhomogeneous Poisson process.

**1. Introduction.** A frequently encountered problem in nonparametric statistics is to estimate a monotone function $\lambda$ on a compact interval, say, $[0,1]$. Grenander [5], Brunk [2] and Huang and Wellner [9] propose estimators defined as the slope of a concave (or convex) approximation of an estimator of a primitive of $\lambda$, in the cases where $\lambda$ is a monotone density function, a monotone regression mean and a monotone failure rate, respectively. These estimators have aroused great interest since they are nonparametric, data driven (they do not require the choice of a smoothing parameter) and easy to implement using, for example, the pool adjacent violators algorithm; see [1]. Moreover, Reboul [14] provides nonasymptotic control of their $\mathbb{L}_1$-risk, which proves that they are optimal in some sense. From an asymptotic point of view, Prakasa Rao [13], Brunk [3] and Huang and Wellner [9] prove cube-root convergence of these estimators at a fixed point and obtain the pointwise asymptotic distribution; Groeneboom, Hooghiemstra and Lopuhaä [8] and Durot [4] prove a central limit theorem for the $\mathbb{L}_1$-error of the Grenander and Brunk estimators, respectively, and Kulikov and Lopuhaä [11] generalize the result in [8] to the $\mathbb{L}_p$-error of the Grenander estimator.

In this paper we consider the problem of estimating a monotone function $\lambda : [0,1] \to \mathbb{R}$ in a general model. We provide a unified approach for studying the $\mathbb{L}_p$-error of estimators defined as the slope of a concave (or convex) approximation of an estimator of a primitive of $\lambda$. We prove that, at a point that may depend on the number $n$ of observations and is far enough from 0 and 1, the local $\mathbb{L}_p$-risk is of order $n^{-p/3}$. We also provide control of the local $\mathbb{L}_p$-risk near the boundaries and derive the result that the global $\mathbb{L}_p$-risk is of order $n^{-p/3}$. Our main result is a central limit theorem for the $\mathbb{L}_p$-error; see Theorem 2: we prove that the $\mathbb{L}_p$-error is asymptotically Gaussian with explicit (though unknown) asymptotic mean and variance. Applying the results to the regression and density models, we recover the results of [4, 8, 11] about Brunk and Grenander estimators. Also, we obtain new results for the Huang–Wellner estimator in the random censorship model, and for an estimator of a monotone intensity function based on $n$ independent copies of an inhomogeneous Poisson process. We believe that our method applies to other models.

Our main motivation for proving asymptotic normality of the $\mathbb{L}_p$-error relies on goodness-of-fit tests. Assume indeed we wish to test $H_0 : \lambda = \lambda_0$ for a given decreasing (resp. increasing) $\lambda_0$, against the nonparametric alternative that $\lambda$ is decreasing (resp. increasing). Using asymptotic normality and proper estimators for the asymptotic mean and variance, we can draw from the observations a normalization of the $\mathbb{L}_p$-distance between $\hat{\lambda}_n$ and $\lambda_0$ that converges under $H_0$ to the standard Gaussian law. The test that rejects $H_0$ if this normalization exceeds the $(1-\alpha)$-quantile of the standard Gaussian law has asymptotic level $\alpha$. With additional effort, Theorem 2 can also be used to test a composite null hypothesis. This will be detailed elsewhere.

The paper is organized as follows. In Section 2 we define and study our estimator in a general model. In Section 3 we apply the results of Section 2 to the random censorship, inhomogeneous Poisson process, regression and density models. The results of Section 2 are proved in Sections 4 and 5 and the results of Section 3 are proved in Section 6.

**2. Main results.** We aim at estimating a function $\lambda : [0,1] \to \mathbb{R}$ subject to the constraint that it is nonincreasing (or nondecreasing), on the basis of $n$ observations. Assume we have at hand a cadlag (i.e., right continuous with left-hand limits at every point) step estimator $\Lambda_n$ of

$$\Lambda(t) = \int_0^t \lambda(u)\, du, \qquad t \in [0,1].$$

We define the monotone estimator of $\lambda$ as follows:

DEFINITION 1. Let $\Lambda_n : [0,1] \to \mathbb{R}$ be a cadlag step process. If $\lambda$ is nonincreasing (resp. nondecreasing), then the monotone estimator $\hat{\lambda}_n$ based on $\Lambda_n$ is defined as the left-hand slope of the least concave majorant (resp. greatest convex minorant) of $\Lambda_n$, with $\hat{\lambda}_n(0) = \lim_{t \downarrow 0} \hat{\lambda}_n(t)$.

Thus, the monotone estimator is a step process that can jump only at the jump points of $\Lambda_n$; it is monotone and left-continuous.

Hereafter, $M_n$ denotes the process defined on $[0,1]$ by $M_n = \Lambda_n - \Lambda$. We make the following assumptions.

(A1) $\lambda$ is monotone and differentiable on $[0,1]$ with $\inf_t |\lambda'(t)| > 0$ and $\sup_t |\lambda'(t)| < \infty$.

(A2) There exists $C > 0$ such that, for all $x \geq n^{-1/3}$ and $t \in [0,1]$,

$$\mathbb{E}\left[\sup_{u \in [0,1],\ x/2 \leq |t-u| \leq x} (M_n(u) - M_n(t))^2\right] \leq \frac{Cx}{n}. \tag{1}$$

(A2′) Inequality (1) holds for all $x > 0$ and $t \in \{0,1\}$.

First, we give a control of the local $\mathbb{L}_p$-risk of $\hat{\lambda}_n$ at a time $t$ that is allowed to depend on $n$: it is of order $n^{-p/3}$ if $t$ is far enough from 0 and 1 [in particular, if $t \in (0,1)$ does not depend on $n$]. We obtain a control of larger order if $t$ is near a boundary and derive a control of the global $\mathbb{L}_p$-risk:

THEOREM 1. *Assume* (A1), (A2), (A2′) *and let* $p \in [1,2)$. *Then there exists* $K > 0$, *which depends only on* $\lambda$, $C$ *and* $p$, *such that*

$$\mathbb{E}|\hat{\lambda}_n(t) - \lambda(t)|^p \leq K n^{-p/3}$$

*for all* $t \in [n^{-1/3}, 1 - n^{-1/3}]$, *and*

$$\mathbb{E}|\hat{\lambda}_n(t) - \lambda(t)|^p \leq K[n(t \wedge (1-t))]^{-p/2} \tag{2}$$

*for all* $t \in (0, n^{-1/3}] \cup [1 - n^{-1/3}, 1)$.

COROLLARY 1. *Assume* (A1), (A2), (A2′) *and let* $p \in [1,2)$. *Then*

$$\mathbb{E}\left[\int_0^1 |\hat{\lambda}_n(t) - \lambda(t)|^p\, dt\right] = O(n^{-p/3}).$$

Note that Theorem 1 does not provide a control of the risk at $t \in \{0,1\}$. In fact, it is known that the monotone estimator is not consistent at the points 0 and 1 in particular models; see [17] for the density model. To control the error at the boundaries, we assume the following.

(A3) $\hat{\lambda}_n(0)$ and $\hat{\lambda}_n(1)$ are stochastically bounded.

The following lemma provides a sufficient condition for (A3), which will be useful for applications.

LEMMA 1. *Assume (A1), (A2) and (A2′). If for every $\varepsilon > 0$ there exists $\delta > 0$ such that the probability that $\Lambda_n$ jumps in $(0, \delta/n)$ or in $(1 - \delta/n, 1)$ is less than $\varepsilon$, then (A3) holds.*

PROOF. Let $x$, $\delta$ and $\varepsilon$ be fixed positive numbers. One has

$$\mathbb{P}(|\hat{\lambda}_n(0)| > x) \leq \mathbb{P}(|\hat{\lambda}_n(\delta/n)| > x) + \mathbb{P}(\hat{\lambda}_n(0) \neq \hat{\lambda}_n(\delta/n)).$$

From Theorem 1, $\hat{\lambda}_n(\delta/n)$ is stochastically bounded. Moreover, $\hat{\lambda}_n(0)$ can differ from $\hat{\lambda}_n(\delta/n)$ only if $\Lambda_n$ jumps in $(0, \delta/n)$. Hence, both probabilities in the above upper bound are less than $\varepsilon$, provided $\delta$ is small enough and $x$ is large enough, whence $\hat{\lambda}_n(0) = O_\mathbb{P}(1)$. Likewise, $\hat{\lambda}_n(1) = O_\mathbb{P}(1)$. □

To compute the asymptotic distribution of the $\mathbb{L}_p$-error, we assume that $M_n$ can be approximated in distribution by a Gaussian process. Specifically, we assume the following.

(A4) Let $B_n$ be either a Brownian bridge or a Brownian motion. There exist $q > 12$, $C_q > 0$, $L : [0, 1] \to \mathbb{R}$ and versions of $M_n$ and $B_n$ such that

$$\mathbb{P}\left(n^{1 - 1/q} \sup_{t \in [0,1]} |M_n(t) - n^{-1/2} B_n \circ L(t)| > x\right) \leq C_q x^{-q}$$

for all $x \in (0, n]$. Moreover, $L$ is increasing and twice differentiable on $[0, 1]$ with $\sup_t |L''(t)| < \infty$ and $\inf_t L'(t) > 0$.

We also need to define the following process $X$:

(3) $$X(a) = \arg\max_{u \in \mathbb{R}}\{-(u - a)^2 + W(u)\}, \qquad a \in \mathbb{R},$$

where $W$ is a standard two-sided Brownian motion (see [6, 7] for a precise description of this process). It is known that, for every $p > 0$, $\mathbb{E}|X(0)|^p$ is finite and the following number $k_p$ is well defined and finite:

$$k_p = \int_0^\infty \text{cov}(|X(0)|^p, |X(a) - a|^p)\, da.$$

We are now in position to state our main result.

THEOREM 2. *Assume (A1), (A2′), (A3) and (A4). Assume, moreover, there are $C' > 0$ and $s > 3/4$ with*

(4) $$|\lambda'(t) - \lambda'(x)| \leq C'|t - x|^s \qquad \text{for all } t, x \in [0, 1].$$

Let $p \in [1, 5/2)$. Then with $m_p = \mathbb{E}|X(0)|^p \int_0^1 |4\lambda'(t)L'(t)|^{p/3} \, dt$,

$$n^{1/6}\left(n^{p/3}\int_0^1 |\hat{\lambda}_n(t) - \lambda(t)|^p \, dt - m_p\right)$$

converges in distribution as $n \to \infty$ to the Gaussian law with mean zero and variance $\sigma_p^2 = 8k_p \int_0^1 |4\lambda'(t)L'(t)|^{2(p-1)/3} L'(t) \, dt$.

Note that our proof of Theorem 2 is partly inspired by [4, 8, 11]. As in those papers, a key step consists in proving that the $\mathbb{L}_p$-error of $\hat{\lambda}_n$ is asymptotically equivalent to an $\mathbb{L}_p$-error of $\hat{U}_n$, the inverse process of $\hat{\lambda}_n$. In the present approach, the proof is quite simple (even for $p > 1$) thanks to the use of Theorem 1. Another key step consists in approximating a proper normalization of $\hat{U}_n(a)$ by the location of the maximum of a drifted Brownian motion. In the present approach, thanks to Proposition 1 in [4], we deal with a parabolic drift independent of $n$, whereas in [8] and [11] the considered drift depends on $n$ and is only close to parabolic (which brings about technicalities, e.g., in the computation of asymptotic moments). Finally, asymptotic normality is proved using Bernstein's method of big blocks and small blocks, as in [8] and [11].

Let us comment on the assumptions in Theorem 2. On one hand, the contribution of the boundaries of the $\mathbb{L}_p$-error is not negligible for $p \geq 5/2$ because $\hat{\lambda}_n$ converges slowly to $\lambda$ near 0 and 1 (this was already stressed for the density model in [11]). This is the reason why we restrict ourself to $p < 5/2$. On the other hand, our proof of Theorem 2 relies on Proposition 1 of [4], which provides a control of the error we make when we approximate the location of the maximum of a given process by that of a drifted Brownian motion. The assumptions $q > 12$ and $s > 3/4$ emerge when using this proposition; see Lemma 5 below. We believe that the proposition can be improved with the assumptions $q > 12$ and $s > 3/4$ being weakened.

To conclude this section, we comment on a slight modification of $\hat{\lambda}_n$. Let $\mathcal{C}_n$ be the set consisting of 0, 1 and the jump points of $\Lambda_n$, and let $C_n$ be the "cumulative sum diagram" consisting of the points $(t, \Lambda_n(t))$, $t \in \mathcal{C}_n$. If $\lambda$ is nonincreasing (resp. nondecreasing), let $\tilde{\lambda}_n$ be the left-hand slope of the least concave majorant (resp. greatest convex minorant) of $C_n$. Then $\hat{\lambda}_n$ and $\tilde{\lambda}_n$ are identical if $\Lambda_n$ is nondecreasing and $\lambda$ is nonincreasing, but they may differ otherwise. In some applications, $\tilde{\lambda}_n$ may be preferred to $\hat{\lambda}_n$ since, for instance, the least-squares estimator of a monotone regression mean takes the form $\tilde{\lambda}_n$. Therefore, we now describe the asymptotic behavior of $\tilde{\lambda}_n$. Let $\tilde{\Lambda}_n$ be the continuous piecewise-affine version of $\Lambda_n$, which means that $\tilde{\Lambda}_n(t) = \Lambda_n(t)$ at every $t \in \mathcal{C}_n$, and $\Lambda_n$ is affine in between two consecutive such points. Assume

(5) $$\sup_{t \in [0,1]} \mathbb{E}[(\tilde{\Lambda}_n(t) - \Lambda_n(t))^2] \leq Cn^{-4/3}$$

for some $C > 0$. Then Theorem 1 and Corollary 1 remain true with $\hat{\lambda}_n$ replaced by $\tilde{\lambda}_n$. On the other hand, assume

$$
(6) \qquad \mathbb{E}\left[\sup_{t \in [0,1]} |\tilde{\Lambda}_n(t) - \Lambda_n(t)|^q\right] \leq C n^{1-q}
$$

for some $q > 12$ and $C > 0$. Then Theorem 2 remains true with $\hat{\lambda}_n$ replaced by $\tilde{\lambda}_n$. The proof of these results is omitted. It is worth noticing that the extra assumptions (5) and (6) hold in every application we consider in Section 3.

**3. Applications.** In this section we consider several models where it may be interesting to estimate a function $\lambda$ on $[0,1]$ subject to a monotonicity constraint. In each model we propose an estimator $\Lambda_n$ of $\Lambda$, we give sufficient conditions for the assumptions (A2), (A2′), (A3) and (A4), and we make explicit the function $L$ in (A4). In particular, this provides sufficient conditions for the $\mathbb{L}_p$-error of the monotone estimator to be asymptotically Gaussian with explicit asymptotic mean and variance. It is worth noticing that, in each considered application, (A2) and (A2′) follow from Doob's inequality and the fact that a proper modification of $M_n$ is a martingale. Also, (A4) follows from an embedding argument similar to that of Komlós, Major and Tusnády [10].

3.1. *The random censorship model.* Assume we observe a right-censored sample $(X_1, \delta_1), \ldots, (X_n, \delta_n)$. Here, $X_i = \min(T_i, Y_i)$ and $\delta_i = \mathbb{1}_{T_i \leq Y_i}$, where the $T_i$'s are nonnegative i.i.d. failure times and the $Y_i$'s are i.i.d. censoring times independent of the $T_i$'s. Assume that the common distribution function $F$ of the $T_i$'s is absolutely continuous with density function $f$ and that we aim at estimating the failure rate $\lambda = f/(1-F)$ on $[0,1]$. Let $N_n$ be the Nelson–Aalen estimator, defined as follows: if $t_1 < \cdots < t_k$ are the distinct times when we observe uncensored data and $n_i$ is the number of $X_j$ that are greater than or equal to $t_i$, then $N_n$ is constant on each $[t_i, t_{i+1})$ with

$$
N_n(t_i) = \sum_{j \leq i} \frac{1}{n_j}.
$$

Moreover, $N_n(t) = 0$ for all $t < t_1$ and $N_n(t) = N_n(t_k)$ for all $t \geq t_k$. Let $\Lambda_n$ be the restriction of $N_n$ to $[0,1]$ and $G$ be the common distribution function of the $Y_i$'s. The monotone estimator based on $\Lambda_n$ is the Huang–Wellner estimator and we have the following.

THEOREM 3. *Assume* (A1), $F(1) < 1$ *and* $\lim_{t \uparrow 1} G(t) < 1$.

  (i) *Then* (A2), (A2′) *and* (A3) *hold.*

(ii) *Assume, moreover,* $\inf_{t\in[0,1]} \lambda(t) > 0$ *and $G$ has a bounded continuous first derivative on $(0,1)$. Then* (A4) *holds with*

$$L(t) = \int_0^t \frac{\lambda(u)}{(1-F(u))(1-G(u))} \, du, \qquad t \in [0,1]. \tag{7}$$

Note that in the case of nonrandom censoring times $Y_i \equiv 1$, one has $G(u) = 0$ for all $u < 1$, so $L$ reduces to $L = (1-F)^{-1} - 1$.

### 3.2. The Poisson process model.
Assume we observe i.i.d. inhomogeneous Poisson processes $N_1, \ldots, N_n$, and their common mean function $\Lambda$ is differentiable on $[0,1]$ with derivative $\lambda$. Let $\Lambda_n$ be the restriction of $\sum_i N_i/n$ to $[0,1]$. Then we have the following.

THEOREM 4. *Assume* (A1), $\Lambda(1) < \infty$ *and* $\inf_{t\in[0,1]} \lambda(t) > 0$. *Then* (A2), (A2'), (A3) *and* (A4) *hold with* $L = \Lambda$.

### 3.3. The regression model.
Assume we observe $y_{i,n} = \lambda(i/n) + \varepsilon_{i,n}$, $i = 1, \ldots, n$, where the $\varepsilon_{i,n}$'s are independent random variables with mean zero. Let

$$\Lambda_n(t) = \frac{1}{n} \sum_{i \leq nt} y_{i,n}, \qquad t \in [0,1].$$

Then the monotone estimator based on $\Lambda_n$ is (a slight modification of) the Brunk estimator and we have the following.

THEOREM 5. *Assume* (A1) *and* $\sup_{i,n} \mathbb{E}|\varepsilon_{i,n}|^q \leq c_q$ *for some* $q \geq 2$ *and* $c_q > 0$.

(i) *Then* (A2), (A2') *and* (A3) *hold.*
(ii) *Assume, moreover,* $q > 12$ *and* $\mathrm{var}(\varepsilon_{i,n}) = \sigma^2(i/n)$ *for some* $\sigma^2 : [0,1] \to \mathbb{R}_+$. *If $\sigma^2$ has a bounded first derivative and satisfies* $\inf_t \sigma^2(t) > 0$, *then* (A4) *holds with* $L(t) = \int_0^t \sigma^2(u) \, du$.

In particular, if the $\varepsilon_{i,n}$'s are i.i.d. with a finite moment of order $q > 12$ and variance $\sigma^2 > 0$, then $L$ reduces to $L(t) = t\sigma^2$. Thus, we recover Theorems 1 and 2 of [4].

### 3.4. The density model.
Assume we observe independent random variables $X_1, \ldots, X_n \in [0,1]$ with common distribution function $\Lambda$ and density function $\lambda = \Lambda'$. Then, the monotone estimator based on the empirical distribution function of $X_1, \ldots, X_n$ is the Grenander estimator and we have the following.

THEOREM 6. *Assume* (A1) *and* $\inf_t \lambda(t) > 0$. *Then* (A2), (A2′), (A3) *and* (A4) *hold with* $L = \Lambda$.

In particular, we recover Theorem 1.1 of [8] and Theorem 1.1 of [11].

**4. Proof of Theorem 1.** We assume here $\lambda$ is decreasing. The similar proof in the increasing case is omitted. We denote by $K$ or $K'$ (resp. $c$) a positive number that depends only on $\lambda$, $C$ and $p$ and that can be chosen as large (resp. small) as we wish. The same letter may denote different constants in different formulas.

First, we give upper bounds for the tail probabilities of the inverse process. Recall that for every nonincreasing left-continuous function $h:[0,1] \to \mathbb{R}$, the (generalized) inverse of $h$ is defined as follows: for every $a \in \mathbb{R}$, $h^{-1}(a)$ is the greatest $t \in [0,1]$ that satisfies $h(t) \geq a$, with the convention that the supremum of an empty set is zero. Let $\Lambda_n^+$ be the upper version of $\Lambda_n$ defined as follows: $\Lambda_n^+(0) = \Lambda_n(0)$ and for every $t \in (0,1]$,

$$\Lambda_n^+(t) = \max\left\{\Lambda_n(t), \lim_{u \uparrow t} \Lambda_n(u)\right\}.$$

Setting $\hat{U}_n = (\hat{\lambda}_n)^{-1}$, one can check that

(8) $\qquad \hat{U}_n(a) = \underset{u \in [0,1]}{\arg\max}\{\Lambda_n^+(u) - au\} \qquad$ for all $a \in \mathbb{R}$,

where arg max denotes the greatest location of the maximum (which is achieved). Moreover, for any $a \in \mathbb{R}$ and $t \in (0,1]$, one has $\hat{U}_n(a) \geq t$ if and only if $a \leq \hat{\lambda}_n(t)$. Hereafter, $g = \lambda^{-1}$.

LEMMA 2. *There exists* $K > 0$ *such that, for every* $a \in \mathbb{R}$ *and* $x > 0$,

(9) $\qquad \mathbb{P}[|\hat{U}_n(a) - g(a)| \geq x] \leq \dfrac{K}{nx^3}.$

PROOF. Fix $a \in \mathbb{R}$, $x \geq n^{-1/3}$ and denote by $P_x$ the probability in (9). By (8), we can have $|\hat{U}_n(a) - g(a)| > x/2$ only if there exists $u \in [0,1]$ with $|u - g(a)| > x/2$ and $\Lambda_n^+(u) - au \geq \Lambda_n^+(g(a)) - ag(a)$, whence

$$P_x \leq \mathbb{P}\left[\sup_{|u-g(a)|>x/2}\{\Lambda_n^+(u) - au\} \geq \Lambda_n^+(g(a)) - ag(a)\right].$$

But $\Lambda_n$ is cadlag and $\Lambda_n^+ \geq \Lambda_n$, so the previous inequality remains true with $\Lambda_n^+$ replaced by $\Lambda_n$. Let $c$ satisfy $0 < c < \inf_t |\lambda'(t)|/2$. If $\lambda(g(a)) \neq a$, then either $a > \lambda(g(a))$ and $g(a) = 0$, or $a < \lambda(g(a))$ and $g(a) = 1$. Hence, from Taylor's expansion,

(10) $\qquad \Lambda(u) - \Lambda(g(a)) \leq (u - g(a))a - c(u - g(a))^2$

for all $u \in [0,1]$, whence

(11) $\quad P_x \leq \mathbb{P}\left[\sup_{|u-g(a)|>x/2}\{M_n(u) - M_n(g(a)) - c(u-g(a))^2\} \geq 0\right].$

It then follows from Markov's inequality and (A2) that

$$P_x \leq \sum_{k \geq 0} \mathbb{P}\left[\sup_{|u-g(a)|\in[x2^{k-1},x2^k]}\{M_n(u) - M_n(g(a))\} \geq c(x2^{k-1})^2\right]$$

$$\leq C \sum_{k \geq 0} \frac{x2^k/n}{(cx^2 2^{2k-2})^2}.$$

But $\sum_k 2^{-3k}$ is finite so (9) holds for all $x \geq n^{-1/3}$. This inequality clearly extends to all $x > 0$ since the upper bound is greater than one for all $x < n^{-1/3}$, provided $K \geq 1$. □

LEMMA 3. *There exists $K > 0$ such that, for every $x > 0$ and $a \notin \lambda([0,1])$,*

(12) $\quad \mathbb{P}[|\hat{U}_n(a) - g(a)| \geq x] \leq \frac{K}{nx(\lambda(g(a))-a)^2}.$

PROOF. We argue as above except that we use (A2′) instead of (A2), and instead of (10), we use the fact that $\Lambda(u) - \Lambda(g(a)) \leq (u-g(a))\lambda(g(a))$. □

Now we prove Theorem 1. Let $t \in (0,1)$. By the Fubini theorem,

$$I_1 := \mathbb{E}[(\hat{\lambda}_n(t) - \lambda(t))_+]^p = \int_0^\infty \mathbb{P}[\hat{\lambda}_n(t) - \lambda(t) > x]px^{p-1}\,dx,$$

where for all $x \in \mathbb{R}$, $x_+ = \max(x,0)$. We have $\hat{U}_n(\lambda(t)+x) \geq t$ whenever $\hat{\lambda}_n(t) > \lambda(t) + x$, whence

$$I_1 \leq \int_0^\infty \mathbb{P}[\hat{U}_n(\lambda(t)+x) \geq t]px^{p-1}\,dx.$$

By (A1), there exists $c > 0$ such that $g(\lambda(t)+x) \leq t - cx$ for every number $x$ that satisfies $\lambda(t) + x \in (\lambda(1), \lambda(0))$. As a probability is no more than one, it thus follows from Lemma 2 that

(13) $\quad I_1 \leq Kn^{-p/3} + \int_{\lambda(0)-\lambda(t)}^\infty \mathbb{P}[\hat{U}_n(\lambda(t)+x) \geq t]px^{p-1}\,dx.$

One has $g(\lambda(t)+x) = 0$ for all $x > \lambda(0) - \lambda(t)$, so Lemma 2 yields

$$\mathbb{P}[\hat{U}_n(\lambda(t)+x) \geq t] \leq \frac{K}{nt^3}.$$

Assume $t \geq n^{-1/3}$. Combining this with (12) yields

$$I_1 \leq Kn^{-p/3} + \int_{\lambda(0)-\lambda(t)+n^{-1/3}}^{\infty} \frac{K}{nt(\lambda(0)-\lambda(t)-x)^2} px^{p-1}\,dx.$$

As $p < 2$, we obtain $I_1 \leq Kn^{-p/3}$. Now assume $t \leq n^{-1/3}$. Then $n^{-1/3} \leq (nt)^{-1/2}$. A probability is no more than one, so (13) and (12) yield

$$I_1 \leq K(nt)^{-p/2} + \int_{\lambda(0)-\lambda(t)+(nt)^{-1/2}}^{\infty} \frac{K}{nt(\lambda(0)-\lambda(t)-x)^2} px^{p-1}\,dx.$$

As $p < 2$, we obtain $I_1 \leq K(nt)^{-p/2}$. In both cases,

$$I_1 \leq K(n^{-p/3} + (nt)^{-p/2}).$$

Similarly,

$$I_2 := \mathbb{E}[(\lambda(t) - \hat{\lambda}_n(t))_+]^p \leq K(n^{-p/3} + (n(1-t))^{-p/2})$$

and the result follows.

**5. Proof of Theorem 2.** We assume here $\lambda$ is decreasing. The similar proof in the increasing case is omitted. We denote by $K$ or $K'$ (resp. $c$) a positive number that depends only on $\lambda$, $C$, $p$, $C_q$, $q$, $L$, and that can be chosen as large (resp. small) as we wish. The same letter may denote different constants in different formulas. Moreover, we denote by $\hat{U}_n$ the inverse process (8) and we set $g = \lambda^{-1}$. We first provide in Lemma 4 an upper bound for the tail probability of $\hat{U}_n$, which is sharper than (9). Then, thanks to Proposition 1 in [4], we prove two lemmas that will be useful to approximate a properly normalized version of $\hat{U}_n(a)$ with the location of the maximum of a drifted Brownian motion. Finally, we prove Theorem 2.

LEMMA 4. *There exists $K > 0$ such that, for every $a \in \mathbb{R}$ and $x > 0$,*

(14) $$\mathbb{P}[|\hat{U}_n(a) - g(a)| \geq x] \leq K(nx^3)^{1-q}.$$

PROOF. Fix $a \in \mathbb{R}$, $x \in (0,1]$ and denote by $P_x$ the probability in (14). From (11), one has $P_x \leq P'_x + P''_x$, where $P'_x$ is equal to

$$\mathbb{P}\left(\sup_{|u-g(a)|>x/2} \left\{n^{-1/2}(B_n \circ L(u) - B_n \circ L(g(a))) - \frac{c}{2}(u-g(a))^2\right\} \geq 0\right)$$

and

$$P''_x = \mathbb{P}\left(\sup_{u \in [0,1]} |M_n(u) - n^{-1/2} B_n \circ L(u)| \geq \frac{cx^2}{16}\right).$$

One can derive from the properties of Brownian motion and the Brownian bridge (see, e.g., (24) below and the proof of Theorem 4 in [4]) that, for all $x \in (0,1]$,

$$P'_x \leq K \exp(-cnx^3) \leq K'(nx^3)^{1-q}.$$

Now by (A4), there exists $K > 0$ with

$$P''_x \leq K x^{-2q} n^{1-q} \leq K(nx^3)^{1-q}.$$

Hence, (14) holds for all $x \in (0,1]$. It clearly extends to all $x > 0$ since both $\hat{U}_n(a)$ and $g(a)$ belong to $[0,1]$. □

LEMMA 5. *Let $T_n > 0$, $W_n$ be a standard two-sided Brownian motion, $D_n : [-T_n, T_n] \to \mathbb{R}$ a nonrandom function and $R_n$ a process indexed by $[-T_n, T_n]$. Furthermore, let*

$$U_n = \arg\max_{[-T_n, T_n]} \{D_n + W_n + R_n\} \quad \text{and} \quad V_n = \arg\max_{[-\log n, \log n]} \{D_n + W_n\}.$$

*Assume $D_n$ continuously differentiable, $D_n(0) = 0$ and there exist positive $A$ and $c$ such that $|D'_n(u)| \leq A|u|$ and $D_n(u) \leq -cu^2$ for all $u \in [-T_n, T_n]$. Assume, moreover, either (i) or (ii), where:*

(i) $T_n = n^{1/(3(6q-11))}$ *for some $q > 12$ and there exists $K > 0$ such that*

(15) $$\mathbb{P}\left[\sup_{u \in [-T_n, T_n]} |R_n(u)| > x\right] \leq K x^{-q} n^{1-q/3} \quad \text{for all } x \in (0, n^{2/3}].$$

(ii) $T_n = \log n$ *and there exist $K > 0$ and $s > 3/4$ with*

$$\sup_{u \in [-T_n, T_n]} |R_n(u)| \leq K n^{-s/3} (\log n)^3.$$

*Let $r = 2(q-1)/(2q-3)$ under (i) and $r < 2s$ under (ii). Then there exists $K' > 0$ that depends only on $K$, $A$, $c$ and $r$ such that*

$$\mathbb{E}|U_n - V_n|^r \leq K' \left(\frac{n^{-1/6}}{\log n}\right)^r.$$

PROOF. Assume (i). Assume, moreover, $n$ is large enough so that $T_n \geq \log n$. If $V'_n$ denotes the greatest location of the maximum of $D_n + W_n$ on $[-T_n, T_n]$, then $V_n$ can differ from $V'_n$ only if $|V'_n| > \log n$. It thus follows from Proposition 1 in [4] (see also the comments just above this proposition) that there exists an absolute constant $C$ such that the probability that $|U_n - V_n| > \delta$ is no more than

$$\mathbb{P}\left[\sup_{u \in [-T_n, T_n]} |R_n(u)| > \frac{x \delta^{3/2}}{2}\right] + Cx \log n + 2\mathbb{P}(|V'_n| > \log n)$$

for every $(x, \delta)$ that satisfies

(16) $\quad \delta \in (0, \log n], \ x > 0, \quad A^2 (\log n)^2 \leq \dfrac{1}{2\delta \log(1/2x\delta)}.$

Moreover, for every $x \geq 0$,

(17) $\quad \mathbb{P}(|V_n'| \geq x) \leq 2\exp(-c^2 x^3/2);$

see, for example, Theorem 4 in [4]. Let $\varepsilon > 0$ and for every $\delta > 0$, set

$$x_\delta = (\log n)^{-1/(q+1)} \delta^{-3q/(2(q+1))} n^{(3-q)/(3(q+1))}.$$

Then (16) holds for every $(\delta, x_\delta)$ with $\delta \in (n^{-1/6}/\log n, n^{-\varepsilon})$, provided $n$ is large enough. By (15), there thus exists $K' > 0$ such that, for every such $\delta$,

(18) $\quad \mathbb{P}(|U_n - V_n| > \delta) \leq K' x_\delta \log n.$

Now, $|U_n - V_n| \leq 2T_n$, so from Fubini's theorem

$$\mathbb{E}|U_n - V_n|^r = \int_0^{2T_n} \mathbb{P}(|U_n - V_n| > \delta) r \delta^{r-1} \, d\delta.$$

But for every $\delta > n^{-\varepsilon}$, $|U_n - V_n|$ can exceed $\delta$ only if it exceeds $n^{-\varepsilon}$ and therefore, the above integral is no more than

$$\left(\dfrac{n^{-1/6}}{\log n}\right)^r + K' \log n \int_{n^{-1/6}/\log n}^{n^{-\varepsilon}} x_\delta r \delta^{r-1} \, d\delta + K' x_{n^{-\varepsilon}} \log n (2T_n)^r.$$

Since $r < 3q/(2(q+1))$, straightforward computations prove that this is of order $O((n^{-1/6}/\log n)^r)$, provided $q > 12$ and $\varepsilon$ is small enough. This completes the proof in the case (i).

Assume (ii). For every $\delta > 0$, let

$$x_\delta = 2K \delta^{-3/2} n^{-s/3} (\log n)^3.$$

Arguing as above, we get (18) for every $\delta \in (n^{-1/6}/\log n, n^{-\varepsilon})$. We conclude with the same arguments, since $s > 3/4$ and $r < 2s$. □

LEMMA 6. *Let $U_n$ and $V_n$ be processes indexed by $J_n \subset [x_0, x_1]$ for some real numbers $x_0$ and $x_1$ independent of $n$. Let $p \geq 1$, $r > 1$ and let $r'$ satisfy $1/r = 1 - 1/r'$. Assume there are $q'$ and $K$ such that*

(19) $\quad \sup_{a \in J_n} \mathbb{E}|U_n(a)|^{q'} \leq K \quad and \quad \sup_{a \in J_n} \mathbb{E}|V_n(a)|^{q'} \leq K$

*for all $n$. Assume, moreover, either* (i) *or* (ii), *where:*

(i) $q' = (p-1)r'$ *and* $\sup_{a \in J_n} \mathbb{E}|U_n(a) - V_n(a)|^r = o(n^{-r/6})$.
(ii) $q' = pr'$ *and there exist $\gamma > r/6$ and $K' > 0$ such that, for every $n$ and $a \in J_n$, $\mathbb{P}(U_n(a) \neq V_n(a)) \leq K' n^{-\gamma}$.*

*Then*

$$\int_{J_n} |U_n(a)|^p \, da = \int_{J_n} |V_n(a)|^p \, da + o_{\mathbb{P}}(n^{-1/6}).$$

PROOF. It follows from Taylor's expansion that

(20) $$|x^p - y^p| \leq p|x - y|(x \vee y)^{p-1} \leq p|x - y|(x^{p-1} + y^{p-1})$$

for all positive numbers $x$ and $y$. Hence, for every $a \in J_n$,

$$\mathbb{E}||U_n(a)|^p - |V_n(a)|^p| \leq p\mathbb{E}[|U_n(a) - V_n(a)|(|U_n(a)|^{p-1} + |V_n(a)|^{p-1})].$$

Also,

$$\mathbb{E}||U_n(a)|^p - |V_n(a)|^p| \leq \mathbb{E}[\mathbb{1}_{U_n(a) \neq V_n(a)}(|U_n(a)|^p + |V_n(a)|^p)].$$

Hence, the result follows from Hölder's inequality. □

Now we turn to the proof of the theorem. Hereafter,

$$\mathcal{J}_n = n^{p/3} \int_0^1 |\hat{\lambda}_n(t) - \lambda(t)|^p \, dt.$$

• *Step* 1. First we express $\mathcal{J}_n$ in terms of $\hat{U}_n$. Precisely, we prove

(21) $$\mathcal{J}_n = n^{p/3} \int_{\lambda(1)}^{\lambda(0)} |\hat{U}_n(a) - g(a)|^p |g'(a)|^{1-p} \, da + o_{\mathbb{P}}(n^{-1/6}).$$

For every $x \in \mathbb{R}$, let $x_+ = \max(x, 0)$. Moreover, let

$$I_1 = \int_0^1 [(\hat{\lambda}_n(t) - \lambda(t))_+]^p \, dt, \qquad I_2 = \int_0^1 [(\lambda(t) - \hat{\lambda}_n(t))_+]^p \, dt$$

and

$$J_1 = \int_0^1 \int_0^{(\lambda(0) - \lambda(t))^p} \mathbb{1}_{\hat{\lambda}_n(t) \geq \lambda(t) + a^{1/p}} \, da \, dt.$$

We have $\hat{\lambda}_n(t) < \lambda(0)$ for all $t > \hat{U}_n(\lambda(0))$, so

$$0 \leq I_1 - J_1 = \int_0^{\hat{U}_n(\lambda(0))} \int_{(\lambda(0) - \lambda(t))^p}^{\infty} \mathbb{1}_{\hat{\lambda}_n(t) \geq \lambda(t) + a^{1/p}} \, da \, dt$$

$$\leq \int_0^{\hat{U}_n(\lambda(0))} [(\hat{\lambda}_n(t) - \lambda(t))_+]^p \, dt.$$

Hence, by monotonicity

$$I_1 - J_1 \leq \int_0^{n^{-1/3} \log n} (\hat{\lambda}_n(t) - \lambda(t))_+^p \, dt + |\hat{\lambda}_n(0) - \lambda(1)|^p \mathbb{1}_{n^{1/3} \hat{U}_n(\lambda(0)) > \log n}.$$

Let $p' \in (p-1/2, 2)$ be such that $1 \leq p' \leq p$ (such a $p'$ exists since $p \in [1, 5/2)$). By assumption, $\lambda$ is bounded and $\hat{\lambda}_n(0)$ is stochastically bounded, so, from Lemma 4,

$$I_1 - J_1 \leq |\hat{\lambda}_n(0) - \lambda(1)|^{p-p'} \int_0^{n^{-1/3}\log n} |\hat{\lambda}_n(t) - \lambda(t)|^{p'} dt + o_\mathbb{P}(n^{-p/3-1/6}).$$

Now, note that the results in Theorem 1 remain true under the assumptions of Theorem 2 (since one can use Lemma 4 instead of Lemma 2 in the proof). As $p' \in [1, 2)$, we get

$$\mathbb{E}\left(\int_0^{n^{-1/3}\log n} |\hat{\lambda}_n(t) - \lambda(t)|^{p'} dt\right) \leq K n^{-(1+p')/3} \log n.$$

But $p' > p - 1/2$, so

$$\int_0^{n^{-1/3}\log n} |\hat{\lambda}_n(t) - \lambda(t)|^{p'} dt = o_\mathbb{P}(n^{-p/3-1/6}).$$

Therefore, $I_1 = J_1 + o_\mathbb{P}(n^{-p/3-1/6})$. The change of variable $b = \lambda(t) + a^{1/p}$ then yields

$$I_1 = \int_{\lambda(1)}^{\lambda(0)} \int_{g(b)}^{\hat{U}_n(b)} p(b - \lambda(t))^{p-1} \mathbb{1}_{g(b) < \hat{U}_n(b)} \, dt \, db + o_\mathbb{P}(n^{-p/3-1/6}).$$

By Taylor's expansion, (A1) and (4), there exists $K > 0$ such that

(22) $\quad |[b - \lambda(t)]^{p-1} - [(g(b) - t)\lambda' \circ g(b)]^{p-1}| \leq K(t - g(b))^{p-1+s}$

for all $b \in (\lambda(1), \lambda(0))$ and $t \in (g(b), 1)$. As a probability is no more than one, integrating (14) proves that, for every $q' < 3(q-1)$, there exists $K_{q'} > 0$ with

(23) $\quad \mathbb{E}[(n^{1/3}|\hat{U}_n(a) - g(a)|)^{q'}] \leq K_{q'} \qquad$ for all $a \in \mathbb{R}$.

Thus, $I_1$ equals

$$\int_{\lambda(1)}^{\lambda(0)} \int_{g(b)}^{\hat{U}_n(b)} p(t - g(b))^{p-1} |\lambda' \circ g(b)|^{p-1} \mathbb{1}_{g(b) < \hat{U}_n(b)} \, dt \, db + R + o_\mathbb{P}(n^{-p/3-1/6}),$$

where $R = O_\mathbb{P}(n^{-(p+s)/3})$. Hence,

$$I_1 = \int_{\lambda(1)}^{\lambda(0)} |\hat{U}_n(b) - g(b)|^p |\lambda' \circ g(b)|^{p-1} \mathbb{1}_{g(b) < \hat{U}_n(b)} \, db + o_\mathbb{P}(n^{-p/3-1/6}).$$

Likewise,

$$I_2 = \int_{\lambda(1)}^{\lambda(0)} |g(b) - \hat{U}_n(b)|^p |\lambda' \circ g(b)|^{p-1} \mathbb{1}_{g(b) > \hat{U}_n(b)} \, db + o_\mathbb{P}(n^{-p/3-1/6})$$

and the result follows, since $\mathcal{J}_n = n^{p/3}(I_1 + I_2)$.

• *Step* 2. Now we approximate a proper normalization of $\hat{U}_n$ by $\tilde{V}$, defined as follows. We have the representation

$$B_n(t) = W_n(t) - \xi_n t, \tag{24}$$

where $W_n$ is a standard Brownian motion, $\xi_n \equiv 0$ if $B_n$ is a Brownian motion and $\xi_n$ is a standard Gaussian variable independent of $B_n$ if $B_n$ is a Brownian bridge. Let $d = |\lambda'|/2(L')^2$, and for every $t \in [0,1]$ let

$$W_t(u) = n^{1/6}[W_n(L(t) + n^{-1/3}u) - W_n(L(t))], \tag{25}$$

so that $W_t$ is a standard Brownian motion. For every $t \in [0,1]$, we define $\tilde{V}(t)$ as the location of the maximum of the drifted Brownian motion $u \mapsto -d(t)u^2 + W_t(u)$ over $[-\log n, \log n]$. We aim at proving

$$\mathcal{J}_n = \int_0^1 \left| \tilde{V}(t) - n^{-1/6} \frac{\xi_n}{2d(t)} \right|^p \left| \frac{\lambda'(t)}{L'(t)} \right|^p dt + o_{\mathbb{P}}(n^{-1/6}). \tag{26}$$

For every $a \in \mathbb{R}$, let $a^\xi = a - n^{-1/2}\xi_n L'(g(a))$. The process $\hat{U}_n$ is nonincreasing and $|\xi_n|$ is less than $\log n$ with probability greater than $1 - \exp(-(\log n)^2/2)$. As $L'$ is bounded, we derive from Lemma 4 that

$$\mathbb{P}(|L(\hat{U}_n(a^\xi)) - L(g(a))| > x) \le K(nx^3)^{1-q} \tag{27}$$

for all $x \in [n^{-1/3}, L(1) - L(0)]$ and large enough $n$. With a modification of $K$, this inequality holds for all $x > 0$ and $n \in \mathbb{N}$. As a probability is no more than one, integrating this inequality yields

$$\sup_{a \in \mathbb{R}} \mathbb{E}[(n^{1/3}|L(\hat{U}_n(a^\xi)) - L(g(a))|)^{q'}] \le K, \tag{28}$$

provided $q' < 3(q-1)$. Recall (21). Then Lemma 6(i) with, for example, $r = r' = 2$ combined with Hölder's inequality and the change of variable $a \to a^\xi$ proves that

$$\mathcal{J}_n = n^{p/3} \int_{\lambda(1)}^{\lambda(0)} \left| \frac{L(\hat{U}_n(a)) - L(g(a))}{L'(g(a))} \right|^p |g'(a)|^{1-p} \, da + o_{\mathbb{P}}(n^{-1/6})$$

$$= n^{p/3} \int_{J_n} |L(\hat{U}_n(a^\xi)) - L(g(a^\xi))|^p \frac{|g'(a)|^{1-p}}{(L'(g(a)))^p} \, da + o_{\mathbb{P}}(n^{-1/6}),$$

where

$$J_n = [\lambda(1) + n^{-1/6}/\log n, \lambda(0) - n^{-1/6}/\log n].$$

Let $a \in \mathbb{R}$. By (8),

$$L(\hat{U}_n(a^\xi)) = \arg\max_{u \in [L(0), L(1)]} \{(\Lambda_n^+ \circ L^{-1} - a^\xi L^{-1})(u)\}.$$

The location of the maximum of a process $\{Z(u),\ u \in I\}$ is also the location of the maximum of $\{AZ(u)+B,\ u \in I\}$ for any $A>0$ and $B \in \mathbb{R}$. Therefore,

$$n^{1/3}(L(\hat{U}_n(a^\xi)) - L(g(a))) = \arg\max_{u \in I_n(a)}\{D_n(a,u) + W_{g(a)}(u) + R_n(a,u)\},$$

where $W_{g(a)}$ is given by (25),

$$I_n(a) = [-n^{1/3}(L(g(a)) - L(0))\ n^{1/3}(L(1) - L(g(a)))],$$

$$D_n(a,u) = n^{2/3}(\Lambda \circ L^{-1} - aL^{-1})(L(g(a)) + n^{-1/3}u) - n^{2/3}(\Lambda(g(a)) - ag(a))$$

and $R_n(a,u)$ is equal to

$$n^{2/3}(a - a^\xi)(L^{-1}(L(g(a)) + n^{-1/3}u) - g(a)) - n^{-1/6}\xi_n u + \tilde{R}_n(a,u)$$

for some $\tilde{R}_n$ which satisfies

$$\sup_{a \in \mathbb{R}, u \in I_n(a)} |\tilde{R}_n(a,u)| \le n^{2/3} \sup_{t \in [0,1]} |\Lambda_n^+(t) - \Lambda(t) - n^{-1/2}B_n \circ L(t)|.$$

We will use Lemma 5 to show that $R_n$ is negligible. For this task, we need to localize. Let $T_n = n^{1/(3(6q-11))}$ and

$$\tilde{U}_n(a) = \arg\max_{u \in [-T_n, T_n]} \{D_n(a,u) + W_{g(a)}(u) + R_n(a,u)\}.$$

If $n$ is large enough, then $[-T_n, T_n] \subset I_n(a)$ for all $a \in J_n$, so

$$n^{1/3}(L(\hat{U}_n(a^\xi)) - L(g(a)))$$

can differ from $\tilde{U}_n(a)$ only if its absolute value exceeds $T_n$. It thus follows from (27) and (28) that we can apply Lemma 6(ii) with some $r' < 3(q-1)/p$, $r'$ as close as possible to $3(q-1)/p$. We get

$$\mathcal{J}_n = \int_{J_n} |\tilde{U}_n(a) + n^{1/3}(L(g(a)) - L(g(a^\xi)))|^p \frac{|g'(a)|^{1-p}}{(L'(g(a)))^p}\,da + o_\mathbb{P}(n^{-1/6})$$

$$= \int_{J_n} \left|\tilde{U}_n(a) - n^{-1/6}\frac{\xi_n}{2d(g(a))}\right|^p \frac{|g'(a)|^{1-p}}{(L'(g(a)))^p}\,da + o_\mathbb{P}(n^{-1/6}).$$

Now let

$$\tilde{\tilde{U}}_n(a) = \arg\max_{u \in [-\log n, \log n]} \{D_n(a,u) + W_{g(a)}(u)\}.$$

By Taylor's expansion, there are positive $K$ and $c$ with

$$\left|\frac{\partial}{\partial u}D_n(a,u)\right| \le K|u| \quad \text{and} \quad D_n(a,u) \le -cu^2$$

for every $a \in J_n$ and $u \in [-T_n, T_n]$. Moreover, there exists $K > 0$ with

$$|R_n(a,u)| \le Ku^2 n^{-1/2}|\xi_n| + n^{2/3}\sup_{t \in [0,1]}|\Lambda_n(t) - \Lambda(t) - n^{-1/2}B_n \circ L(t)|,$$

since $\Lambda_n$ is cadlag. By (A4), (15) thus holds with $R_n(u)$ replaced by $R_n(a,u)$. Due to Theorem 4 in [4], $\tilde{\tilde{U}}_n(a)$ has bounded moments of any order, so we can apply Lemmas 5 and 6 both with condition (i) to get

$$\mathcal{J}_n = \int_{J_n} \left| \tilde{\tilde{U}}_n(a) - n^{-1/6} \frac{\xi_n}{2d(g(a))} \right|^p \frac{|g'(a)|^{1-p}}{(L'(g(a)))^p} \, da + o_{\mathbb{P}}(n^{-1/6}).$$

Now we approximate $\tilde{\tilde{U}}_n(a)$ by $\tilde{V}(g(a))$. By Taylor's expansion and (4), there exists $K$ such that, for all $|u| \leq \log n$,

$$|D_n(a,u) - d(g(a))u^2| \leq K n^{-s/3} (\log n)^3.$$

It follows from (17) that $\tilde{V}(t)$ has bounded moments of any order so Lemma 5(ii) and Lemma 6(i) show that

$$\mathcal{J}_n = \int_{J_n} \left| \tilde{V}(g(a)) - n^{-1/6} \frac{\xi_n}{2d(g(a))} \right|^p \frac{|g'(a)|^{1-p}}{(L'(g(a)))^p} \, da + o_{\mathbb{P}}(n^{-1/6}),$$

and (26) follows from the change of variable $t = g(a)$.

• *Step* 3. Now we prove that, although $B_n$ could be a Brownian bridge in (A4), everything works as if it were a Brownian motion. This is similar to Corollary 3.3 in [8] and Lemma 2.2 in [11], but the present argument takes a simpler form since we deal with $\tilde{V}$. Precisely, we show that $\xi_n$ can be removed from (26), that is,

(29) $$\mathcal{J}_n = \int_0^1 |\tilde{V}(t)|^p \left| \frac{\lambda'(t)}{L'(t)} \right|^p dt + o_{\mathbb{P}}(n^{-1/6}).$$

This is precisely (26) if $B_n$ is a Brownian motion since, in that case, $\xi_n \equiv 0$. Hence, we assume here that $B_n$ is a Brownian bridge. Therefore, $\xi_n$ is a standard Gaussian variable. Let

$$\mathcal{D}_n = n^{1/6} \left\{ \int_0^1 |\tilde{V}(t)|^p \left| \frac{\lambda'(t)}{L'(t)} \right|^p dt - \int_0^1 \left| \tilde{V}(t) - n^{-1/6} \frac{\xi_n}{2d(t)} \right|^p \left| \frac{\lambda'(t)}{L'(t)} \right|^p dt \right\}.$$

We will show that $\mathcal{D}_n = o_{\mathbb{P}}(1)$. Hereafter, for every $t$, $V(t)$ denotes the location of the maximum of the process $u \mapsto -d(t)u^2 + W_t(u)$ over $\mathbb{R}$. Then for every $t$, $V(t)$ can differ from $\tilde{V}(t)$ only if $|V(t)| > \log n$, so similar to (17),

(30) $$\mathbb{P}(\tilde{V}(t) \neq V(t)) \leq 2 \exp(-c^2 (\log n)^3).$$

Moreover,

$$d(t)^{2/3} V(t) = \arg\max_{u \in \mathbb{R}} \{-u^2 d(t)^{-1/3} + W_t(u d(t)^{-2/3})\},$$

which, by scaling, is distributed as $X(0)$; see (3). Fix $\gamma \in (0, 1/12)$. Corollaries 3.4 and 3.3 in [7] show that $X(0)$ has a bounded density function, so from (30),

$$\mathbb{P}(|\tilde{V}(t)| \leq n^{-\gamma}) \leq K n^{-\gamma}.$$

Here, $K$ does not depend on $t$ since $d$ is bounded. Moreover, $\xi_n$ and $\tilde{V}(t)$ possess uniformly bounded moments of any order and the probability that $|\xi_n|$ exceeds $\log n$ is less than $\exp(-(\log n)^2/2)$. Expanding $x \mapsto x^p$ around $|\tilde{V}(t)|$ then proves that $\mathcal{D}_n$ is asymptotically equivalent to

$$pn^{1/6}\int_0^1 \left\{|\tilde{V}(t)| - \left|\tilde{V}(t) - n^{-1/6}\frac{\xi_n}{2d(t)}\right|\right\}|\tilde{V}(t)|^{p-1}\left|\frac{\lambda'(t)}{L'(t)}\right|^p \mathbb{1}_{A_n(t)}\,dt,$$

where $A_n(t)$ is the intersection of the events $\{|\tilde{V}(t)| > n^{-\gamma}\}$ and $\{|\xi_n| \leq \log n\}$. Hence,

$$\mathcal{D}_n = p\xi_n \int_0^1 \frac{\tilde{V}(t)}{2d(t)}|\tilde{V}(t)|^{p-2}\left|\frac{\lambda'(t)}{L'(t)}\right|^p dt + o_{\mathbb{P}}(1).$$

Now, $\tilde{V}(t)$ has a symmetric distribution, so

$$\mathbb{E}\left(\int_0^1 \frac{\tilde{V}(t)}{2d(t)}|\tilde{V}(t)|^{p-2}\left|\frac{\lambda'(t)}{L'(t)}\right|^p dt\right)^2 = \mathrm{var}\left(\int_0^1 \frac{\tilde{V}(t)}{2d(t)}|\tilde{V}(t)|^{p-2}\left|\frac{\lambda'(t)}{L'(t)}\right|^p dt\right)$$

and one can prove, arguing as in Step 5 below, that this tends to zero as $n \to \infty$. Thus, the above integral converges to zero in probability. As $\xi_n$ is stochastically bounded, we get $\mathcal{D}_n = o_{\mathbb{P}}(1)$.

• *Step* 4. Now, we prove that it is sufficient to show

$$n^{1/6}\int_0^1 Y_n(t)\,dt \to \mathcal{N}(0, \sigma_p^2) \qquad \text{in distribution,}$$

where

$$Y_n(t) = (|\tilde{V}(t)|^p - \mathbb{E}|\tilde{V}(t)|^p)\left|\frac{\lambda'(t)}{L'(t)}\right|^p.$$

We have seen that $d(t)^{2/3}V(t)$ is distributed as $X(0)$, so (30) implies

$$\int_0^1 \mathbb{E}|\tilde{V}(t)|^p \left|\frac{\lambda'(t)}{L'(t)}\right|^p dt = \mathbb{E}|X(0)|^p \int_0^1 d(t)^{-2p/3}\left|\frac{\lambda'(t)}{L'(t)}\right|^p dt + o(n^{-1/6})$$

$$= m_p + o(n^{-1/6}).$$

Thus, by (29),

$$n^{1/6}(\mathcal{J}_n - m_p) = n^{1/6}\int_0^1 Y_n(t)\,dt + o_{\mathbb{P}}(1),$$

which proves the stated result.

• *Step* 5. In this step we show

(31) $$\lim_{n\to\infty}\mathrm{var}\left(n^{1/6}\int_0^1 Y_n(t)\,dt\right) = \sigma_p^2.$$

Let $v_n = \text{var}(\int_0^1 Y_n(t)\,dt)$. From Fubini's theorem,
$$v_n = 2\int_0^1 \int_s^1 \left|\frac{\lambda'(t)}{L'(t)} \times \frac{\lambda'(s)}{L'(s)}\right|^p \text{cov}(|\tilde{V}(t)|^p, |\tilde{V}(s)|^p)\,dt\,ds.$$

Let $c_n = 2n^{-1/3}\log n/\inf_t L'(t)$. The increments of $W_n$ are independent, so $\tilde{V}(t)$ and $\tilde{V}(s)$ are independent for all $|t-s| \geq c_n$. Moreover, $|\tilde{V}(t)|$ possesses bounded moments of any order, so

(32) $$v_n = 2\int_0^1 \int_s^{\min(1,s+c_n)} \left|\frac{\lambda'(s)}{L'(s)}\right|^{2p} \text{cov}(|\tilde{V}(t)|^p, |\tilde{V}(s)|^p)\,dt\,ds + o(n^{-1/3}).$$

For every $s$ and $t$, let $\tilde{V}_t(s)$ be the location of the maximum of the process $u \mapsto -d(s)u^2 + W_t(u)$ over $[-\log n, \log n]$ and let $V_t(s)$ be the location of the maximum of this process over the whole real line. By (17), $V_t(s)$ and $\tilde{V}_t(s)$ have bounded moments of any order. Hölder's inequality combined with (20) thus yields
$$|\text{cov}(|\tilde{V}_t(t)|^p, |\tilde{V}_s(s)|^p) - \text{cov}(|\tilde{V}_t(s)|^p, |\tilde{V}_s(s)|^p)| \leq K\mathbb{E}^{1/r}|\tilde{V}_t(t) - \tilde{V}_t(s)|^r,$$
where $r > 1$ is arbitrary. Since $\tilde{V}_t(t) = \tilde{V}(t)$, Lemma 5(ii) yields
$$v_n = 2\int_0^1 \int_s^{\min(1,s+c_n)} \left|\frac{\lambda'(s)}{L'(s)}\right|^{2p} \text{cov}(|\tilde{V}_t(s)|^p, |\tilde{V}_s(s)|^p)\,dt\,ds + o(n^{-1/3}).$$

For every fixed $s$, $V_t(s)$ can differ from $\tilde{V}_t(s)$ only if $|V_t(s)| > \log n$, so similar to (17), we get
$$\mathbb{P}(\tilde{V}_t(s) \neq V_t(s)) \leq 2\exp(-c^2(\log n)^3).$$

Thus, $\tilde{V}_t(s)$ and $\tilde{V}_s(s)$ can be replaced by $V_t(s)$ and $V_s(s)$ in the above integral. Now, fix $s$ and $t$ in $[0,1]$ and let $X$ be given by (3), where
$$W(u) = n^{1/6}d(s)^{1/3}(W_n(L(s) + n^{-1/3}d(s)^{-2/3}u) - W_n(L(s))).$$
Then
$$d(s)^{2/3}V_t(s) = X(n^{1/3}d(s)^{2/3}(L(t) - L(s))) - n^{1/3}d(s)^{2/3}(L(t) - L(s)).$$

The change of variable $a = n^{1/3}d(s)^{2/3}(L(t) - L(s))$ and straightforward computations then yield (31).

• *Step* 6. It remains to prove asymptotic normality of $n^{1/6}\int_0^1 Y_n(t)\,dt$. We will use Bernstein's method of big blocks and small blocks, as in [8] and [11]. Let $L_n = n^{-1/3}(\log n)^5$, $L'_n = n^{-1/3}(\log n)^2$ and denote by $N_n$ the integer part of $(L_n + L'_n)^{-1}$. Let $a_0 = 0$, $a_{2N_n+1} = 1$ and for all $n \in \mathbb{N}$ and all $j \in \{0, \ldots, N_n - 1\}$, let $a_{2j+1} = a_{2j} + L_n$ and $a_{2j+2} = a_{2j+1} + L'_n$. Finally, let $\xi_{n,j} = n^{1/6}\int_{a_{2j}}^{a_{2j+1}} Y_n(t)\,dt$. By definition, $\mathbb{E}Y_n(t) = 0$, so
$$\mathbb{E}\left(\sum_{j=0}^{N_n-1} \int_{a_{2j+1}}^{a_{2j+2}} Y_n(t)\,dt\right)^2 = \sum_{i,j} \int_{a_{2j+1}}^{a_{2j+2}} \int_{a_{2i+1}}^{a_{2i+2}} \text{cov}(Y_n(t), Y_n(s))\,dt\,ds.$$

By independence, the terms with $i \neq j$ are equal to zero for large enough $n$, so the above expectation is of order $o(n^{-1/3})$. Hence, $n^{1/6} \int_0^1 Y_n(t)\,dt$ is asymptotically equivalent to $\sum_j \xi_{n,j}$, and by Step 5, $\mathrm{var}(\sum_j \xi_{n,j})$ tends to $\sigma_p^2$ as $n \to \infty$. By Hölder's and Markov's inequalities, we have, for all $\delta > 0$,

$$\sum_{j=0}^{N_n} \mathbb{E}(\xi_{n,j}^2 \mathbb{1}_{|\xi_{n,j}|>\delta}) \leq \sum_{j=0}^{N_n} \mathbb{E}(|\xi_{n,j}|^3)\delta^{-1}.$$

This tends to zero as $n \to \infty$, so the central limit theorem with the Lindeberg condition shows that $\sum_j \xi_{n,j}$ tends to a centered Gaussian distribution with variance $\sigma_p^2$. By Step 4, this completes the proof of the theorem. $\square$

**6. Proof of the results of Section 3.** Here again, $K$, $K'$, $c$, denote positive numbers that do not depend on $n$ and may change from line to line.

6.1. *Proof of Theorem 3.* (i) Let $M_n^*$ be the stopped process

$$M_n^*(t) = M_n(t \wedge X_{(n)}) = \Lambda_n(t) - \Lambda(t \wedge X_{(n)}), \qquad t \in [0,1],$$

where $X_{(n)} = \max_i X_i$. We have $X_{(n)} < 1$ with probability $\gamma^n$, where

(33) $$\gamma = 1 - \lim_{t\uparrow 1}(1 - F(t))(1 - G(t)) < 1.$$

Recall $(a+b)^2 \leq 2a^2 + 2b^2$ for all real numbers $a$ and $b$. As $M_n^*$ is identical to $M_n$ if $X_{(n)} \geq 1$, we get

$$\mathbb{E}\left[\sup_{t \leq u \leq t+x}(M_n(u) - M_n(t))^2\right]$$
$$\leq 2\mathbb{E}\left[\sup_{t \leq u \leq t+x}(M_n^*(u) - M_n^*(t))^2\right] + 2(Kx)^2\gamma^n$$

for every $t \in [0,1]$ and $x \geq 0$. Here, $K$ denotes the supremum norm of $\lambda$. By Theorem 7.5.2 in [16], $M_n^*$ is a square integrable mean zero martingale with predictable variation process

(34) $$\langle M_n^*\rangle(u) = \frac{1}{n}\int_0^u \frac{\lambda(s)}{1 - H_{n-}(s)}\mathbb{1}_{s \leq X_{(n)}}\,ds,$$

where $H_{n-}(s) = n^{-1}\sum_i \mathbb{1}_{X_i < s}$. By Doob's inequality,

$$\mathbb{E}\left[\sup_{t \leq u \leq t+x}(M_n^*(u) - M_n^*(t))^2\right] \leq 4\mathbb{E}[(M_n^*(1 \wedge (t+x)) - M_n^*(t))^2]$$
$$= 4\mathbb{E}[(M_n^*(1 \wedge (t+x)))^2 - (M_n^*(t))^2]$$
$$= \frac{4}{n}\mathbb{E}\left[\int_t^{1\wedge(t+x)} \frac{\lambda(s)}{1 - H_{n-}(s)}\mathbb{1}_{s \leq X_{(n)}}\,ds\right].$$

Let $N$ be the number of $X_i$'s that are greater than or equal to 1. For every $s \leq 1 \wedge X_{(n)}$, $n(1 - H_{n-}(s))$ is greater than or equal to $1 \vee N$. Hence, by monotonicity,

$$\mathbb{E}\left[\sup_{t \leq u \leq t+x} (M_n^*(u) - M_n^*(t))^2\right] \leq 4x\lambda(0)\mathbb{E}\left(\frac{1}{1 \vee N}\right) \leq \frac{Kx}{n},$$

since $N$ has a binomial distribution with parameter $n$ and probability of success $1 - \gamma > 0$. Also, $\gamma^n \leq K/n$ for some $K > 0$, and $x^2 \leq x$ for all $x \in [0,1]$. Hence, for every $t \in [0,1]$ and $x \geq 0$, we have

(35) $$\mathbb{E}\left[\sup_{t \leq u \leq t+x} (M_n(u) - M_n(t))^2\right] \leq \frac{Kx}{n}.$$

To handle the case $u < t$, we derive from (35) that

$$\mathbb{E}\left[\sup_{t-x \leq u \leq t} (M_n(u) - M_n(t))^2\right]$$
$$\leq 2\mathbb{E}[(M_n(t) - M_n((t-x) \vee 0))^2]$$
$$\quad + 2\mathbb{E}\left[\sup_{t-x \leq u \leq t} (M_n(u) - M_n((t-x) \vee 0))^2\right]$$
$$\leq \frac{Kx}{n}$$

for every $t \in [0,1]$ and $x \geq 0$. Combining this with (35) yields (A2) and (A2'). Now, $\Lambda_n$ jumps only at times $t_i$ when we observe uncensored data. Hence, for every $\delta > 0$, the probability that $\Lambda_n$ jumps in $(0, \delta/n)$ or in $(1 - \delta/n, 1)$ is no more than

$$n\mathbb{P}(T_1 \in (0, \delta/n) \cup (1 - \delta/n, 1)).$$

This is no more than $2K\delta$, where $K$ is the supremum norm of $f$ on $[0,1]$, so (A3) follows from Lemma 1.

(ii) Let $L$ be defined by (7), and denote the supremum distance on $[0,1]$ by $\|\cdot\|$. We will prove that there exist versions of $M_n$ and the standard Brownian motion $B_n$ such that, for all $x \in [0, n]$,

(36) $$\mathbb{P}\left[n \sup_{t \in [0,1]} |M_n(t) - n^{-1/2} B_n \circ L(t)| > x + K \log n\right] \leq K' \exp(-cx),$$

where $K$, $K'$ and $c$ depend only on $F$ and $G$. This indeed suffices to prove (A4). We consider the limit-product estimator $F_n$ of Kaplan and Meier,

$$F_n(t) = 1 - \prod_{i \leq k} \left(\frac{n_i - 1}{n_i}\right)^{\mathbb{1}_{t_i \leq t}}, \qquad t \geq 0,$$

and we set $\bar{\Lambda}_n = -\log(1 - F_n)$. By Corollaries 1 and 2 of [12], there are versions of $F_n$ and $B_n$ such that

(37) $\quad \mathbb{P}[n\|F_n - F - n^{-1/2}(1-F)B_n \circ L\| > x + K\log n] \leq K' \exp(-cx)$

for all $x \geq 0$. Here, $K$, $K'$ and $c$ depend only on $F$ and $G$. As $L$ is bounded on $[0,1]$, we have

$$\mathbb{P}[\|B_n \circ L\| \geq x] \leq \exp(-c'x^2)$$

for some $c' > 0$ and all $x \geq 0$. But $F(1) < 1$ and we have (37), so we can assume without loss of generality that $F_n(1) < 1$ and, therefore, $\bar{\Lambda}_n$ is well defined on the whole interval $[0,1]$. As $\Lambda = -\log(1-F)$, expanding $u \mapsto \exp(-u)$ proves that there are positive $c$, $K$ and $K'$, which depend only on $F$ and $G$, such that

$$\mathbb{P}[n\|(\bar{\Lambda}_n - \Lambda)\exp(-\Lambda) - n^{-1/2}(1-F)B_n \circ L\| > x + K\log n] \leq K'e^{-cx},$$

for all $x \in [0, n]$. Hence,

$$\mathbb{P}[n\|\bar{\Lambda}_n - \Lambda - n^{-1/2} B_n \circ L\| > x + K\log n] \leq K' \exp(-cx),$$

and it remains to show that $\bar{\Lambda}_n$ is close enough to $\Lambda_n$. By Taylor's expansion, one has, for all $i$ with $t_i \in [0,1]$,

$$0 \leq \bar{\Lambda}_n(t_i) - \Lambda_n(t_i) \leq \sum_{j \leq i} \frac{1}{2(n_j - 1)^2} \leq \frac{n}{2(N \vee 1 - 1)^2},$$

where we recall that $N$ is the number of $X_i$'s that are greater than or equal to 1. Both $\Lambda_n$ and $\bar{\Lambda}_n$ are constant on the intervals $[t_i, t_{i+1})$. As $N$ is a binomial variable with parameter $n$ and probability of success $1 - \gamma$ [see (33)], one can then derive from Hoeffding's inequality that

$$\mathbb{P}[n\|\bar{\Lambda}_n - \Lambda_n\| > x] \leq K\exp(-cn) \leq K\exp(-cx),$$

for some $K > 0$, $c > 0$ and all $x \in (K', n]$. The result follows.

6.2. *Proof of Theorem 4.* Fix $t \in [0,1]$ and $x > 0$. As $\Lambda_n - \Lambda$ is a martingale, Doob's inequality yields

(38) $\quad \mathbb{E}\left[\sup_{t \leq u \leq t+x} (M_n(u) - M_n(t))^2\right] \leq 4\mathbb{E}((M_n(1 \wedge (t+x)) - M_n(t))^2).$

But $n(\Lambda_n(1 \wedge (t+x)) - \Lambda_n(t))$ has a Poisson distribution with expectation $n(\Lambda(1 \wedge (t+x)) - \Lambda(t))$. Thus, its variance is bounded by $Knx$, where $K$ is the supremum norm of $\lambda$ on $[0,1]$, and (35) holds for all $x > 0$ and $t \in [0,1]$. We can handle the case $u < t$ as in the proof of Theorem 3, whence (A2) and (A2'). Now, $\Lambda_n$ can jump in $(1 - \delta/n, 1)$ only if at least a process $N_i$ jumps in this interval. But the jumps of $N_i$ have height 1, so for every $\delta > 0$,

$$\mathbb{P}(\Lambda_n \text{ jumps in } (1 - \delta/n, 1)) \leq n\mathbb{P}(N_1(1) - N_1(1 - \delta/n) \geq 1).$$

The variable $N_1(1) - N_1(1-\delta/n)$ has a Poisson distribution with expectation $\Lambda(1) - \Lambda(1-\delta/n)$, so by Markov's inequality,

$$\mathbb{P}(\Lambda_n \text{ jumps in } (1-\delta/n, 1)) \leq K\delta.$$

We can proceed likewise to control the probability that $\Lambda_n$ jumps in $(0, \delta/n)$, so (A3) follows from Lemma 1. It remains to prove (A4). For this task, fix $q \geq 2$ and for every $k = 0, \ldots, n$, let $t_k = k/n$. We have

$$(39) \qquad \mathbb{E}|M_n(t_k) - M_n(t_{k-1})|^q \leq Kn^{-q}$$

for all $k \geq 1$ and some $K > 0$. The increments of $M_n$ are independent, so by Theorem 5 in [15], there are versions of $M_n$ and the standard Brownian motion $B_n$ such that

$$\mathbb{E}\left[\max_{1 \leq k \leq n} |M_n(t_k) - n^{-1/2} B_n(\Lambda(t_k))|^q\right] \leq Kn^{1-q}$$

for some $K > 0$. One then obtains, using (39), monotonicity of $\Lambda_n$ and properties of Brownian motion, that there is a $K > 0$ such that

$$\mathbb{E}\left[\sup_{t \in [0,1]} |M_n(t) - n^{-1/2} B_n(\Lambda(t))|^q\right] \leq Kn^{1-q}.$$

This holds for any $q \geq 2$, hence, in particular, for some $q > 12$. Thus, from Markov's inequality, (A4) holds with $L = \Lambda$. $\square$

6.3. *Proof of Theorem* 5. (i) We have $(u+v)^2 \leq 2(u^2 + v^2)$ for all real numbers $u$ and $v$. Hence, for all $t \in [0,1]$ and $x > 0$,

$$\mathbb{E}\left[\sup_{|u-t| \leq x} (M_n(u) - M_n(t))^2\right] \leq \frac{2}{n^2} \mathbb{E}\left[\sup_{|u| \leq x} \left(\sum_{i \leq n(t+u)} \varepsilon_{i,n} - \sum_{i \leq nt} \varepsilon_{i,n}\right)^2\right] + \frac{K}{n^2}$$

for some $K > 0$, which depends only on $\lambda$. By Doob's inequality, this is no more than

$$\frac{8}{n^2} \mathbb{E}\left(\sum_{nt < i \leq n(t+x)} \varepsilon_{i,n}\right)^2 + \frac{8}{n^2} \mathbb{E}\left(\sum_{n(t-x) < i \leq nt} \varepsilon_{i,n}\right)^2 + \frac{K}{n^2} \leq \frac{K'x}{n}$$

for all $x \geq 1/n$, whence (A2). By definition, $\Lambda_n$ jumps at times $i/n$, $i = 1, \ldots, n$. If $t \in \{0, 1\}$, we thus have for every $x \in (0, 1/n)$ that

$$\sup_{|t-u| \leq x} |M_n(u) - M_n(t)| = \sup_{|t-u| \leq x} |\Lambda(u) - \Lambda(t)| \leq Kx,$$

whence (A2′). Moreover, it is clear from Lemma 1 that (A3) holds.

(ii) From Theorem 5 in [15], there exist versions of $(\varepsilon_{i,n})$ and the standard Brownian motion $B_n$ such that

$$\mathbb{E}\left[\sup_{t\in[0,1]}\left|\frac{1}{n}\sum_{i\leq nt}\varepsilon_{i,n} - n^{-1/2}B_n\left(\frac{1}{n}\sum_{i\leq nt}\sigma^2(i/n)\right)\right|^q\right] \leq Kn^{1-q}.$$

Thanks to Markov's inequality, one can then derive (A4) from properties of $B_n$ and the regularity assumptions on $\sigma^2$.

6.4. *Proof of Theorem* 6. Fix $t \in [0,1]$, $x > 0$, and define

$$\mathcal{M}_n(u) = \frac{\Lambda_n(u) - \Lambda_n(t)}{\Lambda(u) - \Lambda(t)}, \qquad u \in [0,1],$$

where we recall that $\Lambda_n$ is the empirical distribution function of the sample $X_1, \ldots, X_n$. By Lemma 2.2 in [8], the process $\{\mathcal{M}_n(u),\ u \in (t,1]\}$ is a reverse time martingale conditionally on $\Lambda_n(t)$. Since $\Lambda$ is increasing and $\lambda$ is bounded, Doob's inequality yields

$$\mathbb{E}\left[\sup_{x/2 \leq u-t \leq x}(M_n(u) - M_n(t))^2\right] \leq Kx^2 \mathbb{E}\left(\frac{M_n(t+x/2) - M_n(t)}{\Lambda(t+x/2) - \Lambda(t)}\right)^2.$$

But $n(\Lambda_n(t+x/2) - \Lambda_n(t))$ is a binomial variable with parameter $n$ and probability of success $\Lambda(t+x/2) - \Lambda(t)$. Moreover, $\lambda$ is bounded away from zero, whence

$$\mathbb{E}\left[\sup_{x/2 \leq u-t \leq x}(M_n(u) - M_n(t))^2\right] \leq \frac{Kx}{n}$$

for all $x > 0$ and $t \in [0,1]$. To handle the case $u < t$, we use the fact that the process $\{\mathcal{M}_n(u),\ u \in [0,t)\}$ is a forward time martingale conditionally on $\Lambda_n(t)$ (see Lemma 2.2 in [8]). Whence, (A2) and (A2′). Now, $\Lambda_n$ jumps at times $X_1, \ldots, X_n$. As $\lambda$ is bounded, the probability that $\Lambda_n$ jumps in $(0, \delta/n)$ or in $(1 - \delta/n, 1)$ is no more than

$$n\mathbb{P}(X_1 \in (0,\delta/n) \cup (1-\delta/n,1))$$

for every $\delta > 0$. This is no more than $2K\delta$, where $K$ is the supremum norm of $\lambda$, so (A3) follows from Lemma 1. Finally, it follows from the Hungarian embedding of [10] that (A4) holds with $L = \Lambda$ and $B_n$ a Brownian bridge.

**Acknowledgment.** I wish to thank Laurence Reboul for useful discussions about the Nelson–Aalen estimator and martingales.

Laboratoire de Statistiques
Université Paris sud
Bâtiment 425
91405 Orsay cedex
France
E-mail: cecile.durot@math.u-psud.fr